\documentclass[a4paper, 12pt]{article}

\usepackage[koi8-r,cp1251]{inputenc}
\usepackage{amsfonts,amssymb,amsmath, hyperref}
\usepackage[final]{epsfig}
\usepackage{graphicx}

\begin{document}

\begin{center}
ON THE DESCRIPTION OF MULTIDIMENSIONAL NORMAL HAUSDORFF
OPERATORS ON LEBESGUE SPACES
\end{center}

\
\begin{center}
A. R. Mirotin
\end{center}

\

\textsc{Abstract.} Hausdorff operators  originated from some classical summation methods. Now this is an active research field. In the present article,  a spectral representation for multidimensional normal Hausdorff operator is given. We show that  normal Hausdorff operator in $L^2(\mathbb{R}^n)$
 is unitary equivalent to the operator of multiplication
by some matrix-valued function (its matrix symbol) in the space $L^2(\mathbb{R}^n; \mathbb{C}^{2^n}).$
 Several  corollaries that show
that properties of a Hausdorff operator are closely related to the properties
of its symbol are considered. In particular,   the norm and the spectrum
of such operators are described in terms of the symbol.

\



\

\section{ Introduction}

\

The notion of a Hausdorff operator with respect to a positive measure on
the unit interval was introduced by Hardy \cite[Chapter XI]{H} as a continuous
analog of the Hausdorff  summability methods for series.  This class of operators contains some important examples, such as Hardy operator,  the Ces\`{a}ro operator and its q-calculus version, and there adjoints. As mentioned in \cite{RF} the Riemann-Liouville fractional integral and the
Hardy-Littlewood-Polya operator   can also be reduced to  the Hausdorff operator, and as was noted in \cite{CDFZ} in the  one-dimensional case the Hausdorff operator is closely related to a Calder\'{o}n-Zygmund convolution operator, too.

The modern theory of
Hausdorff operators begins with  the paper by Liflyand and Moricz \cite{LM}. Now this is an active research field. The survey articles \cite{Ls}, \cite{CFW} contain main results on Hausdorff operators
and bibliography up to 2014. For more resent results see, e.~g., \cite{TAMS}, \cite{CDFZ}, \cite{arx1}, and \cite{JMAA}.
The last paper is devoted to  generalizations of Hausdorff operators to
locally compact groups. In this paper, we accept the following special case
of the definition from \cite{JMAA}.

\

\textbf{Definition 1.} Let  $(\Omega,\mu)$ be some $\sigma$-compact  topological space endowed with a positive regular  Borel measure $\mu,$ $K$  a locally integrable function on $\Omega,$  and $(A(u))_{u\in \Omega}$  a $\mu$-measurable family of $n\times n$-matrices that are nonsingular
for $\mu$-almost every $u$ with $K(u) \ne 0.$
  We define the \textit{Hausdorff  operator} with the kernel $K$  by ($x\in\mathbb{R}^n$ is a column vector)
$$
(\mathcal{H}f)(x)=(\mathcal{H}_{K, A}f)(x) =\int_\Omega K(u)f(A(u)x)d\mu(u).
$$

To our knowledge, all known results on Hausdorff operators refer to the
boundedness of such operators in various settings only (exceptions are the papers  \cite{BHS}, and \cite{R} in which some spectra were calculated for the one-dimensional
case). In particular multidimensional normal Hausdorff operators have not been studied. Our  main goal  is to obtain a spectral representation for such operators. As is known, an explicit diagonalization of a normal operator can be obtained only in a few cases. In this paper,
 using the $n$-dimensional Mellin transform we show that normal Hausdorff operator in $L^2(\mathbb{R}^n)$
with self-adjoint $A(u)$ is unitary equivalent to the operator of multiplication
by some matrix-valued function (its matrix symbol) in the space $L^2(\mathbb{R}^n; \mathbb{C}^{2^n}).$ This
is an analogue of the Spectral Theorem for the class of operators under consideration. This allows us to find the norm and to study  the spectrum
of such operators. The cases of positive definite and negative  definite   $A(u)$ have been examined. We give also for the case of normal operators a negative
answer to the problem of compactness of (nontrivial) Hausdorff operators
posed by Liflyand \cite{L} (see also \cite{Ls}). Several other corollaries are considered that show
that properties of a Hausdorff operator are closely related to the properties
of its symbol. Several examples are worked out, as well.  The results were announced in \cite{arx2}.
It should be noted that the case of $L^p$ spaces is more complicated for the lack of Plancherel Theorem, cf. \cite{arx1}, \cite{Hilb}.

In the following, we assume that all the conditions of definition 1 are fulfilled.

\

\section{ Notation and preliminaries}

   Let $U_j\ (j = 1; \dots ; 2^n)$ be some  fixed enumeration of the family of all open hyperoctants in $\mathbb{R}^n.$  For every
pair $(i; j)$ of indices there is a unique $\varepsilon(i; j)\in\{-1,1\}^n$ such that $\varepsilon(i; j)U_i :=
\{(\varepsilon(i; j)_1 x_1;\dots ; \varepsilon(i; j)_n x_n) : x\in U_i\} = U_j.$ It is clear that $\varepsilon(i; j)U_j = U_i$ and
$\varepsilon(i; j)U_l \cap U_i=\emptyset$ as $l \ne j.$
We will assume that $A(u)$ form a commuting family. Then
as is well known there are an orthogonal $n\times n$-matrix $C$ and a family of
diagonal non-singular real matrices $A'(u) = \mathrm{diag}(a_1(u); \dots ; a_n(u))$ such that
$A'(u) = C^{-1}A(u)C$ for $u\in \Omega.$ Then $a(u) := (a_1(u); \dots ; a_n(u))$ is the family
of all eigenvalues (with their multiplicities) of the matrix $A(u).$ We put
$$
\Omega_{ij} := \{u\in\Omega: (\mathrm{sgn}(a_1(u)); \dots ; \mathrm{sgn}(a_n(u))) = \varepsilon(i; j)\}.
 $$

If $|\det A(u)|^{-1/2}K(u)\in L^1(\Omega)$
we put also
 $$
   \varphi_{ij}(s):=\int_{\Omega_{ij}}K(u)|a(u)|^{-1/2- \imath s}d\mu(u)
   $$
(above we assume that $|a(u)|^{-1/2-\imath s}:=$ $\prod_{j=1}^n |a_j(u)|^{-1/2-\imath s_j}$ where $|a_j(u)|^{-1/2-\imath s_j}:=$ $\exp((-1/2-\imath s_j)\log |a_j(u)|)$).

Evidently, all functions $\varphi_{ij}$ belong to the algebra $C_b(\mathbb{R}^n)$ of bounded and
continuous functions on $\mathbb{R}^n$ and $\varphi_{ij}=\varphi_{ji}.$

\

\textbf{Definition 2.} Let $|a(u)|^{-1/2}K(u)\in L^1(\Omega)$. We define the  \textit{matrix symbol}
of a Hausdorff  operator $\mathcal{H}_{K, A}$  by
$$
\Phi = \left(\varphi_{ij}\right)_{i,j=1}^{2^n}
$$
So, $\Phi$ is a symmetric element of the matrix algebra $\mathrm{Mat}_{2^n}(C_b(\mathbb{R}^n)).$

The symbol was first introduced in \cite{arx1} for the case of positive definite $A(u)$  (in the simplest one-dimensional case the symbol was in fact considered also in \cite[Theorem 2.1]{CDFZ}). As
we shall see properties of a Hausdorff operator are closely related to the
properties of its matrix symbol.

\

\textbf{Remark 1.} Of cause $\Phi$ depends on the enumeration of the family of all  hyperoctants in $\mathbb{R}^n$  we choose.

\

\textbf{Lemma 1} \cite{JMAA} (cf. \cite[(11.18.4)]{H}, \cite{BM}). \textit{Let  $|\det A(u)|^{-1/2}K(u)\in L^1(\Omega).$ Then the operator $\mathcal{H}_{K, A}$ is bounded in $L^2(\mathbb{R}^n)$ and }
$$
\|\mathcal{H}_{K, A}\|\leq \int_\Omega |K(u)||\det A(u)|^{-1/2}d\mu(u).
$$

This estimate is sharp (see \cite[theorem 1]{arx1}).

\

\section{The main result}

\

\textbf{Theorem 1.} \textit{Let $A(u)$ be a commuting family of non-singular self-adjoint
$n\times n$-matrices, and $|\det A(u)|^{-1/2}K(u)\in L^1(\Omega).$ Then the Hausdorff operator
$\mathcal{H}_{K, A}$ in $L^2(\mathbb{R}^n)$ with matrix symbol $\Phi$ is normal and unitary equivalent to the
operator $M_\Phi$ of multiplication by the normal matrix $\Phi$ in the space $L^2(\mathbb{R}^n;\mathbb{C}^{2^n})$ of $\mathbb{C}^{2^n}$-valued
functions. In particular, the spectrum of $\mathcal{H}_{K, A}$ equals to the spectrum $\sigma(\Phi)$
of $\Phi$ in the matrix algebra $\mathrm{Mat}_{2^n}(C_b(\mathbb{R}^n)),$   in other words,}
$$
\sigma(\mathcal{H}_{K, A})=\{\lambda\in \mathbb{C}:\inf_{s\in \mathbb{R}^n}|\det(\lambda-\Phi(s)|=0\}.
$$
\textit{The point spectrum $\sigma_p(\mathcal{H}_{K, A})$ of $\mathcal{H}_{K, A}$ consists of such complex numbers $\lambda$ for which the closed set $E(\lambda):=\{s\in \mathbb{R}^n: \det(\lambda-\Phi(s))=0\}$ has positive Lebesgue measure. The residual spectrum of $\mathcal{H}_{K, A}$  is empty.}

Proof. It is known (see \cite{BM}) that under the conditions of lemma 1 the adjoint
for the Hausdorff operator in $L^2(\mathbb{R}^n)$ has the form
$$
(\mathcal{H}^*f)(x) = (\mathcal{H}_{K, A}^*f)(x) =\int_\Omega \overline{K(v)}|\det A(v)|^{-1}f(A(v)^{-1}x)d\mu(v).
$$
(Thus, the adjoint for a Hausdorff operator is also Hausdorff.) Since $A(u)$
form a commuting family, the normality of $\mathcal{H}_{K, A}$  follows from the equalities
$$
(\mathcal{H}\mathcal{H}^*f)(x)=\int_\Omega \int_\Omega K(u)\overline{K(v)}|\det A(v)|^{-1}f(A(v)^{-1}A(u)x)d\mu(v)d\mu(u),
$$
$$
(\mathcal{H}^*\mathcal{H}f)(x)=\int_\Omega \int_\Omega K(u)\overline{K(v)}|\det A(v)|^{-1}f(A(u)A(v)^{-1}x)d\mu(u)d\mu(v),
$$
and the Fubini Theorem.

Next, let the orthogonal $n\times n$-matrix $C$ and a family of diagonal real
matrices $A'(u) = \mathrm{diag}(a_1(u); \dots ; a_n(u))$ be such that $A'(u)=C^{-1}A(u)C$ for
$u\in \Omega$ as in the section 2. Then all functions $a_j(u)$ are $\mu$-measurable and real,
and $\det A(u)=$ $a_1(u)\dots a_n(u)\ne 0$ for $\mu$-almost all $u.$ Consider the operator
$\hat{C}f(x):=f(Cx)$
in $L^2(\mathbb{R}^n).$ Since
$$
\mathcal{H}_{K, A}=\hat{C}^{-1}\mathcal{H}_{K, A'}\hat{C}
$$
and $\hat{C}$ is unitary, the operator $\mathcal{H}_{K, A}$ is unitary equivalent to $\mathcal{H}_{K, A'}.$

For every pair $(i;j)$ of indices consider the following operator in $L^2(U_j)$:
$$
(H_{ij}f)(x) =\int_{\Omega_{ij}} K(u)f(A'(u)x)d\mu(u).
$$
Then $H_{ij}$ maps $L^2(U_j)$ into $L^2(U_i)$ (because $f(A'(u)x) = 0$ for $f \in L^2(U_j)$
and $x \notin U_i$). Moreover, if $f\in L^2(\mathbb{R}^n)$  and $f_i := f1_{U_i}$ ($1_E$ denotes the indicator of a subset $E\subset \mathbb{R}^n$) then for a.~e. $x \in \mathbb{R}^n$
$$
(\mathcal{H}_{K, A'}f)(x) =\sum\limits_{j=1}^{2^n}\sum\limits_{i=1}^{2^n}(H_{ij}f_i)(x). \eqno(1)
$$
Indeed, for every $x \in \mathbb{R}^n$  that does not belong to any coordinate hyperplane
there is a unique index $j_0$ such that $x \in U_{j_0}.$ Then $A'(u)x\in U_i$ if and only if
$u\in \Omega_{ij_0}.$ Therefore
$$
(\mathcal{H}_{K, A'}f)(x) =\sum\limits_{i=1}^{2^n}\int_\Omega K(u)f_i(A'(u)x)d\mu(u)=\sum\limits_{i=1}^{2^n}\int_{\Omega_{ij_0}}K(u)f_i(A'(u)x)d\mu(u). \eqno(2)
$$
On the other hand, if $x \notin U_j$ then $A'(u)x \notin U_i$ for all $i$ and $u\in
\Omega_{ij}$ and
therefore
$$
\int_{\Omega_{ij}}K(u)f_i(A'(u)x)d\mu(u)=0.
$$
Thus, for $x\in U_{j_0}$ and $f\in L^2(\mathbb{R}^n)$ we have in view of (2) that
$$
\sum\limits_{j=1}^{2^n}\sum\limits_{i=1}^{2^n}(H_{ij}f_i)(x)=
\sum\limits_{j=1}^{2^n}\sum\limits_{i=1}^{2^n}\int_{\Omega_{ij}}K(u)f_i(A'(u)x)d\mu(u)=
$$
$$
\sum\limits_{i=1}^{2^n}\int_{\Omega_{ij_0}}K(u)f_i(A'(u)x)d\mu(u)=(\mathcal{H}_{K, A'}f)(x).
$$
Now, since $\sum_{i=1}^{2^n}H_{ij}f_i\in L^2(U_j)$
for all $j,$ formula (1) can be rewritten as
follows:
$$
\mathcal{H}_{K, A'}f=\bigoplus_{j=1}^{2^n}\sum\limits_{i=1}^{2^n}H_{ij}f_i.\eqno(3)
$$
In turn, if we identify $L^2(\mathbb{R}^n)$ with the orthogonal sum $L^2(U_1)\oplus\dots\oplus L^2(U_{2^n}),$
the last formula can be rewritten in the following way:
$$
\mathcal{H}_{K, A'}f=(H_{ij})_{ij}(f_1;\dots ; f_{2^n})^T
$$
($B^T$ denotes the
transposed to the matrix $B$). So, we get the following block operator matrix representation for $\mathcal{H}_{K, A'}$:
$$
\mathcal{H}_{K, A'}= (H_{ij})_{ij}.
$$

Consider the modified $n$-dimensional Mellin transform for the $n$-hyperoctant
$U_i$ $(i = 1,\dots, 2^n)$:
$$
(\mathcal{M}_if)(s):=\frac{1}{(2\pi)^{n/2}}\int_{U_i}|x|^{-\frac{1}{2}+\imath s}f(x)dx, \quad s\in \mathbb{R}^n.
$$
Then $\mathcal{M}_i$ is a unitary operator acting from $L^2(U_i)$ to $L^2(\mathbb{R}^n).$ This can be easily
obtained from the Plancherel theorem for the Fourier transform by using
an exponential change of variables (see \cite{BPT}). Moreover, if we assume that
$|y|^{-1/2}f(y)\in L^1(U_j)$ then making use of Fubini's theorem, and integrating
by substitution $x=A'(u)^{-1}y=(y_1/a_1(u);\dots;y_n/a_n(u))$ yield the following
($s\in \mathbb{R}^n$):
$$
(\mathcal{M}_i H_{ij}f)(s)=\frac{1}{(2\pi)^{n/2}}\int_{U_i}|x|^{-\frac{1}{2}+\imath s}dx\int_{\Omega_{ij}}K(u)f(A'(u)x)d\mu(u)=
$$
$$
\frac{1}{(2\pi)^{n/2}}\int_{\Omega_{ij}}K(u)d\mu(u)\int_{U_i}|x|^{-\frac{1}{2}+\imath s}f(A'(u)x)dx=
$$
$$
\int_{\Omega_{ij}}K(u)a(u)^{-\frac{1}{2}-\imath s}d\mu(u)\frac{1}{(2\pi)^{n/2}}\int_{U_j}|y|^{-\frac{1}{2}+\imath s}f(y)dy=
\varphi_{ij}(s)(\mathcal{M}_jf)(s).
$$
By continuity we get for all $f\in L^2(U_j)$ that
$$
\mathcal{M}_i H_{ij}f=\varphi_{ij} \mathcal{M}_jf.
$$
So, $H_{ij} =\mathcal{M}_i^{-1} M_{\varphi_{ij}}\mathcal{M}_j$ ($M_{\varphi_{ij}}$ denotes the operator in $L^2(\mathbb{R}^n)$ of multiplication by $\varphi_{ij}$) and therefore
$$
\mathcal{H}_{K, A'}=(\mathcal{M}_i^{-1} M_{\varphi_{ij}}\mathcal{M}_j)_{ij}.
$$

Let $\mathcal{U} := \mathrm{diag}(\mathcal{M}_1; \dots ; \mathcal{M}_{2^n}).$ If we identify $L^2(\mathbb{R}^n)$  with $L^2(U_1)\otimes\dots\otimes L^2(U_{2^n}),$ then $\mathcal{U}$ is a unitary operator between $L^2(\mathbb{R}^n)$  and $L^2(\mathbb{R}^n,\mathbb{C}^{2^n})$  and
$$
\mathcal{H}_{K, A'}=\mathcal{U}^{-1} M_{\Phi}\mathcal{U}.
$$
This proves the first statement of the theorem.

To compute the spectrum, let $\lambda\in \mathbb{C}.$  The operator $\lambda-\mathcal{H}_{K, A}$ is unitary equivalent to the
operator $M_{\lambda-\Phi}$ in $L^2(\mathbb{R}^n,\mathbb{C}^{2^n}).$ The last operator is invertible if and only
if the matrix $\lambda-\Phi(s)$ is invertible (i.~e. $\det(\lambda-\Phi(s))\ne 0$ for all $s \in
\mathbb{R}^n$) and $M_{(\lambda-\Phi)^{-1}}$ acts in $L^2(\mathbb{R}^n,\mathbb{C}^{2^n}).$ This condition is fulfilled if and only
if $(\lambda-\Phi)^{-1}(0;\dots; 0; f; 0; \dots ; 0)^T$ ($f$ is in the $j$th column) belongs to $L^2(\mathbb{R}^n,\mathbb{C}^{2^n})$ for all $f\in L^2(\mathbb{R}^n)$ and
$j = 1, \dots,  2^n.$ Let $(\lambda-\Phi(s))^{-1} = (\alpha_{ij}(\lambda; s))_{ij}.$ Then $\alpha_{ij}(\lambda; \cdot))_{ij} \in C(\mathbb{R}^n)$ and
we have
$$
(\lambda-\Phi)^{-1}(0;\dots; 0; f; 0; \dots ; 0)^T = (\alpha_{ij}(\lambda; s)f)_{i=1}^{2^n}.
$$
So, $\lambda$ is a regular point for $\mathcal{H}_{K, A}$ if and only if the matrix $\lambda-\Phi(s)$ is invertible
and $\alpha_{ij}(\lambda; \cdot) \in C_b(\mathbb{R}^n)$ for every pair $(i;j)$ of indices. This means that the
matrix $\lambda-\Phi$ is invertible in the algebra $\mathrm{Mat}_{2^n}(C_b(\mathbb{R}^n))$ (i.~e. $\lambda\notin \sigma(\Phi)$). But
it is known that the last condition is equivalent to the fact that $\det(\lambda-\Phi)$
is invertible in $C_b(\mathbb{R}^n)$ (see, e.~g., \cite[Prop. VII.3.7, p. 353]{Hu}), i.~e. $\inf_{s\in \mathbb{R}^n}|\det(\lambda-\Phi(s))| > 0.$

Now let $\lambda\in  \sigma_p(M_\Phi)$ and $f\in L^2(\mathbb{R}^n,\mathbb{C}^{2^n})$ a corresponding eigenvalue. Then $(\lambda-\Phi(s))f(s)=0$ for a.e. $s\in \mathbb{R}^n.$ It follows that  $\det(\lambda-\Phi(s))=0$ for a.e. $s\in S$ where the set $S:=\{s\in\mathbb{R}^n: f(s)\ne 0\}$ has a  positive Lebesgue measure. So, $\mathrm{mes}(E(\lambda))>0.$

To prove the converse,  let $\lambda\in \mathbb{C}$ and $\mathrm{mes}(E(\lambda))>0.$ Consider a multifunction $\Gamma(s):=\mathrm{ker}(\lambda-M_{\Phi(s)})\setminus\{0\}$ on $E(\lambda)$ taking values in the set of all subsets of $\mathbb{C}^{2^n}.$ Then $\Gamma(s)\ne\emptyset$ for all $s\in E(\lambda).$ Moreover,
since the map $(s,x)\mapsto (\lambda-\Phi(s))x$ is continuous on  $E(\lambda)\times \mathbb{C}^{2^n}$, the graph $G_\Gamma:=\{(s,x)\in E(\lambda)\times (\mathbb{C}^{2^n}\setminus\{0\}): (\lambda-\Phi(s))x=0\}$ of  $\Gamma$ is a Borel subset of $E(\lambda)\times \mathbb{C}^{2^n}$ (the disjoint union $G_\Gamma \sqcup (E(\lambda)\times\{0\})$ is closed). By the measurable selection theorem (see, e.g., \cite{Levin}) there is a measurable selection $\xi: E(\lambda)\to \mathbb{C}^{2^n},$ $\xi(s)\in \mathrm{ker}(\lambda-M_{\Phi(s)})\setminus\{0\}.$ Let $\chi_C$ be the indicator of a compact subset $C\subset E(\lambda)$ of positive Lebesgue measure. Then the function $f(s):=(\xi(s)/\|\xi(s)\|)\chi_C(s)$ belongs to $L^2(\mathbb{R}^n,\mathbb{C}^{2^n})$ and  is an eigenvalue of $M_\Phi$ which corresponds to $\lambda.$
 Finally, as is well known,  normal operator  has  empty residual spectrum.
This completes the
proof of the theorem.

\

\section{Corollaries and examples}

\

In the following corollaries we assume that the assumptions and notation
of theorem 1 are fulfilled.

\textbf{Corollary 1.} \textit{The operator $\mathcal{H}_{K, A}$ is invertible if and only if $\inf\limits_{s\in \mathbb{R}^n}|\det\Phi(s)| > 0.$ In this case its inverse is unitary equivalent to the operator $M_{\Phi^{-1}}$ in $L^2(\mathbb{R}^n,\mathbb{C}^{2^n}).$}

\textbf{Corollary 2.} \textit{Let $\mathcal{H}_{K_1, A}$ and $\mathcal{H}_{K_2, B}$ be two Hausdorff operators with the
same measure space $(\Omega,\mu)$ such that $(A(u);B(v))$ is a commuting family
of self-adjoint $n\times  n$-matrices that are nonsingular for $\mu$-almost $u$ and $v$
respectively, and $|\det A(u)|^{-1/2}K_1(u), |\det B(v)|^{-1/2}K_2(v)\in L^1(\Omega).$ Then the
product $\mathcal{H}_{K_1, A}\mathcal{H}_{K_2, B}$ is unitary equivalent to the operator $M_{\Phi_1\Phi_2}$  in $L^2(\mathbb{R}^n,\mathbb{C}^{2^n}).$
($\Phi_2$  denotes the matrix symbol of $\mathcal{H}_{K_2, B}$).
}

Proof. First note that the orthogonal matrix $C$ exists such that both
$A'(u) = C^{-1}A(u)C$ and $B'(v) = C^{-1}B(v)C$ are diagonal. Then the proof of
theorem 1 shows that $\mathcal{V}\mathcal{H}_{K_1, A}\mathcal{V}^{-1}=M_{\Phi_1}$  and
 $\mathcal{V}\mathcal{H}_{K_2, B}\mathcal{V}^{-1}=M_{\Phi_2}$   for some unitary
operator $\mathcal{V}$ from $L^2(\mathbb{R}^n)$ to $L^2(\mathbb{R}^n,\mathbb{C}^{2^n})$
 which depends only on $n$ and $C$ and
the result follows.

\

\textbf{Corollary 3.}
$$
\|\mathcal{H}_{K, A}\|=\max\{|\lambda|:\inf_{s\in \mathbb{R}^n}|\det(\lambda-\Phi(s)|=0\}=\sup_{s\in \mathbb{R}^n}\|\Phi(s)\|,
$$
\textit{where $\|\Phi(s)\|$ stands for the  norm of the operator in $\mathbb{C}^{2^n}$ of multiplication by the matrix $\Phi(s).$}

Proof. The first equality follows from  theorem 1 and the normality of $\mathcal{H}_{K, A}$ (the norm
of this operator equals to its spectral radius) and the second one follows from theorem 1 and the equality $\|M_{\Phi}\|=\sup_{s\in \mathbb{R}^n}\|\Phi(s)\|$ (see \cite[Theorem 4.1.1]{Sim} for more general result).

\

\textbf{Corollary 4.} \textit{The matrix symbol of the adjoint operator $\mathcal{H}_{K, A}^*$ is the
adjoint matrix $\Phi^* = (\overline{\varphi_{ij}}).$}

Proof. The adjoint for a Hausdorff operator $\mathcal{H}_{K, A}$  is also Hausdorff \cite{BM}.
By theorem 1 this adjoint is unitary equivalent to the adjoint for the operator
$M_\Phi,$ i.~e. to $M_{\Phi^*}.$

\

\textbf{Corollary 5.} \textit{The Hausdorff operator $\mathcal{H}_{K, A}$ is selfadjoint (positive, unitary) if and only if the matrix $\Phi(s)$ is selfadjoint (respectively, positive definite, unitary) for all $s\in \mathbb{R}^n.$}

Proof. This follows from corollaries 1 and 4.

\

\textbf{Example 1}. (Discrete   Hausdorff operators; see also example 3 below.) Let $\Omega=\mathbb{Z}_+,$ and $\mu$ be a counting measure. Then the definition 1 takes the form ($f\in L^2(\mathbb{R}^n)$)
$$
(\mathcal{H}_{K, A}f)(x) =\sum_{k=0}^\infty K(k)f(A(k)x).
$$
Assume that $\sum_{k=0}^\infty |K(k)||\det A(k)|^{-1/2}<\infty.$ Then $\mathcal{H}_{K, A}$ is bounded on $L^2(\mathbb{R}^n)$ and
$$
\varphi_{ij}(s)=\sum_{k\in \Omega_{ij}} K(k)|\det A(k)|^{-1/2}|a(k)|^{-\imath s}
$$
where $s=(s_l)\in\mathbb{R}^n$ and the principal values of the exponential functions  are considered.
Since this series converges  on $\mathbb{R}^n$ absolutely and uniformly, $\varphi_{ij}$ is uniformly almost periodic. So,
the matrix symbol of $\mathcal{H}_{K, A}$ is a uniformly almost periodic matrix-valued function.

Assume, in addition, that $A(k)=A^k$ where the matrix $A$ is self-adjoint, but not positive definite and $(\lambda_1,\dots,\lambda_n)$ are all eigenvalues of $A$ (with their multiplicities).
Let $\beta:=(\mathrm{sgn}(\lambda_1); \dots ; \mathrm{sgn}(\lambda_n)).$ Then $\beta\in \{-1,1\}^n$ and
$$
(\mathrm{sgn}(a_1(k)); \dots ; \mathrm{sgn}(a_n(k))) =(\mathrm{sgn}\lambda_1^k; \dots ; \mathrm{sgn}\lambda_n^k)=
\begin{cases}
(1,\dots,1) \text {\ if\ } k\in 2\mathbb{Z}_+,\\
\beta \text {\ if \ } k\in 2\mathbb{Z}_+ +1.
\end{cases}
$$
Let us enumerate $n$-hyperoctants $U_j$ in such a way that $U_{2^{n-1}+i}=\beta U_i$ (the coordinate-wise multiplication) for $i = 1, \dots, 2^{n-1}.$ Then $\Omega_{ii}=2\mathbb{Z},$ $\Omega_{ij}=2\mathbb{Z}+1$ if $|j-i|=2^{n-1},$
and $\Omega_{ij}=\emptyset$  otherwise. It follows that
$$
\varphi_{ii}(s)=\varphi_+(s):=\sum_{m=0}^\infty K(2m)|a(2m)|^{-\frac{1}{2}-\imath s}=\sum_{m=0}^\infty \frac{K(2m)}{|\det(A)|^{m}}\prod_{j=1}^n|\lambda_j|^{-\imath 2ms_j}.
$$
Analogously, if $|j-i|=2^{n-1},$
$$
\varphi_{ij}(s)=\varphi_-(s):=\sum_{m=0}^\infty \frac{K(2m+1)}{|\det(A)|^{\frac{2m+1}{2}}}\prod_{j=1}^n|\lambda_j|^{-\imath (2m+1)s_j},
$$
and $\varphi_{ij}=0$ otherwise. So, the matrix symbol is the following block matrix:
$$\Phi=
\begin{pmatrix}
\varphi_+I_{2^{n-1}}&\varphi_- I_{2^{n-1}}\\ \varphi_- I_{2^{n-1}}&\varphi_+I_{2^{n-1}}
\end{pmatrix}
$$
where $I_{2^{n-1}}$ denotes the identity matrix of order
$2^{n-1}.$ Then for every $\lambda\in \mathbb{C}$
$$\lambda-\Phi=
\begin{pmatrix}
(\lambda-\varphi_+)I_{2^{n-1}}&-\varphi_- I_{2^{n-1}}\\ -\varphi_- I_{2^{n-1}}&(\lambda-\varphi_+)I_{2^{n-1}}
\end{pmatrix}
$$
and therefore by the formula of Schur (see, e.~g., \cite[p. 46]{G}),
$$
\det(\lambda-\Phi) = \det((\lambda-\varphi_+)^2  I_{2^{n-1}}-\varphi_-^2 I_{2^{n-1}}) = ((\lambda-\varphi_+ -\varphi_-)(\lambda-\varphi_++\varphi_-))^{2^{n-1}}=
$$
$$
((\lambda-\varphi)(\lambda-\varphi^*))^{2^{n-1}},
$$
 where $\varphi:=\varphi_+ +\varphi_-,$   $\varphi^*:=\varphi_+ -\varphi_-.$ Theorem 1 implies that (we use the boundedness of $\varphi, \varphi^*$)
$$
\sigma(\mathcal{H}_{K, A})=\{\lambda\in \mathbb{C}: \inf_{s\in \mathbb{R}^n}|(\lambda-\varphi(s))(\lambda-\varphi^*(s))|=0\}= \mathrm{cl}(\varphi(\mathbb{R}^n)\cup \varphi^*(\mathbb{R}^n)).
$$
In view of the normality of $\mathcal{H}_{K, A},$ this implies that $\mathcal{H}_{K, A}=\max\{\sup|\varphi|, \sup|\varphi^*|\}.$

Theorem 1 implies also that
$$
\sigma_p(\mathcal{H}_{K, A})=\{\lambda\in \mathbb{C}:\mathrm{mes}(\varphi^{-1}(\{\lambda\})>0\}\cup\{\lambda\in \mathbb{C}:\mathrm{mes}((\varphi^*)^{-1}(\{\lambda\})>0\}=
$$
$$
\{\lambda\in \mathbb{C}:\mathrm{mes}(\varphi^{-1}(\{\lambda\})\cup(\varphi^*)^{-1}(\{\lambda\})) >0\}.
$$
\

As was mentioned above, the problem of compactness of nontrivial Hausdorff operators was posed in \cite{L}. In our case the answer to this question is
negative (the case of positive definite matrices was considered in \cite{JMAA}, \cite{arx1}).

\textbf{Corollary 6.} \textit{The Hausdorff operator $\mathcal{H}_{K, A}$ is noncompact provided it is non-zero.}

Proof. Let $\mathcal{H}_{K, A}$ be compact in $L^2( \mathbb{R}^n)$ and non-zero. We shall use
notation and formulas from the proof of theorem 1. There is $H_{ij}$ that is nonzero, too. Moreover, $H_{ij}$ is compact because it is equal to $P_i\mathcal{H}_{K, A'}|L^2(U_j)$ by
(3) ($P_i$ denotes the orthogonal projection of $L^2( \mathbb{R}^n)$ onto $L^2(U_i)$). It follows that the
operator $\mathcal{M}_iH_{ij}\mathcal{M}_j^{-1}= M_{\varphi_{ij}}$ is non-zero and compact in $L^2( \mathbb{R}^n),$ as well. A
contradiction.

\textbf{Corollary 7.} \cite{arx1}. \textit{Let the matrices $A(u)$ be positive definite. Then
the operator  $\mathcal{H}_{K, A}$ is unitary equivalent to the operator of coordinate-wise
multiplication by a function $\varphi\in C_b(\mathbb{R}^n)$ (the $\mathrm{scalar\ symbol}$) in the space
 $L^2(\mathbb{R}^n,\mathbb{C}^{2^n}).$ In particular,}

(i) \textit{the spectrum, the point spectrum, and the continuous spectrum of $\mathcal{H}_{K, A}$ equal to the spectrum (i. e. to the closure of the range of $\varphi$), to the point
spectrum, and to the continuous spectrum of the operator $M'_\varphi$
of multiplication
by $\varphi$ in $L^2(\mathbb{R}^n)$ respectively, the residual spectrum of $\mathcal{H}_{K, A}$  is empty;}

(ii) $\|\mathcal{H}_{K, A}\|  = \sup|\varphi|.$

Proof. Indeed, if all the matrices $A(u)$ are positive definite then
$$
(\mathrm{sgn}(a_1(u));\dots ;\mathrm{sgn}(a_n(u)) = (1; \dots ; 1).
$$
It follows that $\Omega_{ii} =\Omega$  and
$\Omega_{ij} = \emptyset$ for $i\ne j.$ Therefore
$\varphi_{ii} = \varphi$ where
$$
   \varphi(s):=\int_{\Omega}K(u)|a(u)|^{-1/2- \imath s}d\mu(u) \eqno(4)
  $$
and $\varphi_{ij} = 0$ for $i\ne j.$  So, $\Phi= \mathrm{diag}(\varphi; \dots ; \varphi)$ and the corollary
follows from theorem 1 and corollary 3.

\

\textbf{Example 2.} Consider the \textit{Ces\`{a}ro operator} in $L^2(\mathbb{R}^n)$ (see, e.~g., \cite{K}):
$$
(\mathcal{C}_{\alpha,n}f)(x) = \alpha\int\limits_0^1 \frac{(1-u)^{\alpha-1}}{u}f\left(\frac{x}{u}\right)du\quad  (\alpha > 0;  x \in \mathbb{R}^n).
$$

This is a bounded Hausdorff operator where
$\Omega = [0; 1]$ is endowed with the
Lebesgue measure, $K(u) = (1-u)^{\alpha-1}/u,$ and $A(u) = \mathrm{diag}(1/u; \dots ; 1/u)$
($u \in (0; 1)$) is a positive definite matrix. Its scalar symbol is ($s = (s_j)\in \mathbb{R}^n$)
$$
 \varphi(s)= \alpha\int\limits_0^1 \frac{(1-u)^{\alpha-1}}{u}u^{n/2+\imath \sum_js_j}du=\frac{\Gamma(\alpha+1)\Gamma(n/2+\imath \sum_js_j)}{\Gamma(\alpha+n/2+\imath \sum_js_j)}.
$$

Since the modulus of the function
$$
\gamma(t) :=\frac{\Gamma(\alpha+1)\Gamma(n/2+\imath t)}{\Gamma(\alpha+n/2+\imath t)},\quad t\in \mathbb{R}
$$
attains its maximum at $t = 0$ (this follows, e.~g., from \cite[Section 12.13, Example 1]{WW}), we get by corollary 7 that
$$
\|\mathcal{C}_{\alpha,n}\|=\sup_s|\varphi(s)|=\gamma(0)=\frac{\Gamma(\alpha+1)\Gamma(n/2)}{\Gamma(\alpha+n/2)}.
$$
Moreover, corollary 7 implies that the spectrum of $\mathcal{C}_{\alpha,n}$ is a curve given by
the range of $\gamma(t),$ $t\in \mathbb{R}\cup \{\infty\}.$

\textbf{Corollary 8.} \textit{Let the matrices $A(u)$ be negative definite. Then the matrix
symbol of the operator  $\mathcal{H}_{K, A}$ for some enumeration of $n$-hyperoctants is the following block matrix:
$$\Phi=
\begin{pmatrix}
\mathrm{O}&\varphi I_{2^{n-1}}\\ \varphi I_{2^{n-1}}&\mathrm{O}
\end{pmatrix}              \eqno(5)
$$
where $\varphi$ is given by formula (4).  Moreover, $\sigma(\mathcal{H}_{K, A})= \mathrm{cl}(-\varphi(\mathbb{R}^n)\cup \varphi(\mathbb{R}^n)).$ In particular, $\|\mathcal{H}_{K, A}\|=\sup|\varphi|.$
}

Proof. Let us enumerate $n$-hyperoctants $U_j$ in such a way that $U_{2^{n-1}+i}=-U_i$ for $i = 1, \dots, 2^{n-1}.$ Since
$$
(\mathrm{sgn}(a_1(u));\dots ; \mathrm{sgn}(a_n(u))) = (-1; \dots ;-1),
$$
it follows that $\Omega_{ij} =\Omega$  for $|j -i|=2^{n-1}$ and
$\Omega_{ij} =\emptyset$ otherwise. Therefore
$\varphi_{ij} = \varphi$ for  $|j -i|=2^{n-1}$ and
$\varphi_{ij} =0$ otherwise. Thus, $\Phi$ is given by (5) and
for $\lambda\in \mathbb{C}$ we have
$$\lambda-\Phi=\begin{pmatrix}
\lambda  I_{2^{n-1}}&-\varphi I_{2^{n-1}}\\ -\varphi I_{2^{n-1}}&\lambda  I_{2^{n-1}}
\end{pmatrix}.
$$
Now as in example 1  the formula of Schur implies
$$
\det(\lambda-\Phi) = \det(\lambda^2  I_{2^{n-1}}-\varphi^2 I_{2^{n-1}}) = (\lambda^2-\varphi^2)^{2^{n-1}}
$$
and by theorem 1 we get
$$
\sigma(\mathcal{H}_{K, A})=\{\lambda\in \mathbb{C}: \inf_{s\in \mathbb{R}^n}|\lambda^2-\varphi^2(s)|=0\}= \mathrm{cl}(-\varphi(\mathbb{R}^n)\cup \varphi(\mathbb{R}^n)).
$$

The valuer of the norm follows from this formula and normality of $\mathcal{H}_{K, A}.$

\textbf{Example 3.} Consider the q-calculus version of a Ces\`{a}ro operator (see, e.~g., \cite{Ernst} for the definition of the q-integral)
$$
(C_qf)(x):=\frac{1}{x}\int_0^xf(t)d_qt:=(1-q)\sum_{k=0}^\infty q^kf(q^kx).
$$
Here  $f\in L^2(\mathbb{R}),$ $x\in \mathbb{R},$ and $q$ is real, $0<|q|<1.$
This is a bounded discrete Hausdorff operator in the sense of example 1, where $n=1,$ $K(k)=(1-q)q^k,$ $A=q, a(k)=q^k.$ Two cases are possible.

1) $0<q<1.$ In this case one can apply  corollary 7. By  formula (4) the scalar symbol is
$$
\varphi(s)=(1-q)\sum_{k=0}^\infty(q^k)^{1/2-\imath s}=\frac{1-q}{1-\sqrt{q}q^{-\imath s}}.
$$
Now corollary 7 implies that
$$
\sigma(C_q)=\left\{\frac{1-q}{1-\sqrt{q}z}: z\in \mathbb{C}, |z|=1\right\}=\{\lambda\in \mathbb{C}: |\lambda-1|=\sqrt{q}\}.
$$
It follows that $\|C_q\|=\sqrt{q}.$ Moreover, the operator $C_q$ is invertible, and its inverse $(C_q^{-1}g)(x)=(g(x)-qg(qx))/(1-q)$
is unitary equivalent to the operator of coordinate-wise
multiplication by a function $1/\varphi$ in the space
 $L^2(\mathbb{R},\mathbb{C}^{2}).$

2) $-1<q<0.$ In this case one can apply  corollary 8. Again by  formula (4) the scalar symbol is
$$
\varphi(s)=(1-q)\sum_{k=0}^\infty q^k((-q)^k)^{-1/2-\imath s}=
$$
$$
(1-q)\sum_{k=0}^\infty (-1)^k((-q)^k)^{1/2-\imath s}=\frac{1-q}{1+\sqrt{-q}(-q)^{-\imath s}}.
$$
Since $\varphi(\mathbb{R})=\{\lambda\in \mathbb{C}: |\lambda-1|=\sqrt{-q}\}$ as in the  case 1 above, corollary 8 implies that
$$
\sigma(C_q)=\{\lambda\in \mathbb{C}: |\lambda\pm 1|=\sqrt{-q}\}.
$$
It follows that $\|C_q\|=\sqrt{-q}.$ The operator $C_q$ is invertible, and its inverse $(C_q^{-1}g)(x)=(g(x)-qg(qx))/(1-q)$
is unitary equivalent to the operator $M_{\Phi^{-1}}$ in the space
 $L^2(\mathbb{R},\mathbb{C}^{2}),$ where (see formula (5))
$$\Phi^{-1}=
\begin{pmatrix}
0&1/\varphi \\ 1/\varphi &0
\end{pmatrix}.
$$

\textsc{Department of Mathematics and   Programming  Technologies,
F. Skorina Gomel State University},  Sovietskaya, 104, 246019, Gomel, Belarus

E-mail address: amirotin@yandex.ru

\end{document}